\def\IR{{\Bbb R}} 
 
\def\IS{{\Bbb S}} 
 
\def\IK{{\Bbb K}}

\def\IC{\Bbb C} 
\def\ID{{\Bbb D}}

\def\zbar{{\overline{z}}} 
 
\def\wbar{{\overline{w}}}

\documentclass[10pt]{article} 

\usepackage[psamsfonts]{amssymb} 
\usepackage{amsfonts,amsmath}  

\usepackage{epsfig,multicol} 
\DeclareMathOperator{\loc}{loc}

\newtheorem{theorem}{Theorem}
\newtheorem{lemma}{Lemma}

\title{Extremal mappings of finite distortion and the Radon-Riesz property.}
\author{Gaven Martin \& Cong Yao\thanks{
Work of both authors partially supported by the New Zealand Marsden Fund.
 \newline
Institute for Advanced Study, 
Massey University,  Auckland,
New Zealand.
\newline
email: g.j.martin@massey.ac.nz \&
c.yao@massey.ac.nz
\newline
\newline
{\bf MSC Subject.}  30C62 31A05 49J10  }
}
\date{}
\begin{document}
\maketitle
\begin{abstract}  
We consider Sobolev mappings $f\in W^{1,q}(\Omega,\IC)$, $1<q<\infty$, between planar domains $\Omega\subset \IC$. We analyse the Radon-Riesz property for convex functionals of the form 
\[f\mapsto \int_\Omega \Phi(|Df(z)|,J(z,f)) \; dz \]
 and show that under certain criteria,  which hold in important cases,   weak convergence in $W_{loc}^{1,q}(\Omega)$ of (for instance) a minimising sequence can be improved to strong convergence. This finds important applications in the minimisation problems for mappings of  finite distortion and the $L^p$ and $Exp$\,-Teichm\"uller theories.
\end{abstract}
\section{Introduction}
Recently geometric function theory has developed strong connections with the calculus of variations and planar nonlinear elasticity by 
identifying a facinating interplay between analysis and topology for mappings of finite distortion.  As a simple example if $f$ is a homeomorphism,  then $f\in W^{1,1}_{loc}(\Omega)$ implies $J(z,f)\in L^{1}_{loc}(\Omega)$ and $f$ is differentiable almost everywhere. Other examples have led to the solution of the Nitsche conjecture \cite{AIMO,IKO1},  higher regularity of extremal monotone mappings,  \cite{IKO2}, and other applications in nonlinear elasticity,  see for instance \cite{IO2} and the references therein.  Using the method of $p$-harmonic replacement based on the Rado-Choquet-Knesser theorem,  Iwaniec and Onninen \cite{IO1} have shown (aside from some minor technical issues) that for each $p\geq 1$, given a weakly converging sequence of homeomorphisms $h_j:\Omega\to \Omega'$ in $W^{1,p}(\Omega)$,  $h_j\rightharpoonup h$,  then there exist diffeomorphisms $\tilde{h_j} :\Omega\to \Omega'$ with  $\tilde{h_j}$ converging to $h$ strongly in $W^{1,p}(\Omega)$ and with $\tilde{h_j}-h\in W^{1,p}_0(\Omega)$. So weak convergence of a sequence is replaced by strong convergence of a ``nicer'' sequence.  As an application, the authors show that in some natural planar models of nonlinear elasticity the minimizers of the energy are strong limits of homeomorphisms. 

\medskip

An issue with this result however is that the initial sequence $\{h_j\}_{j=1}^{\infty}$ may carry information that the approximations $\{\tilde{h_j}\}_{j=1}^{\infty}$ do not,  simply because the sequence $\{\tilde{h_j}\}_{j=1}^{\infty}$ may not be a minimising sequence even if $\{h_j\}_{j=1}^{\infty}$ is.  For instance uniform or local uniform bounds on the higher integrability of certain convex combinations of minors may not follow from any $W^{1,p}$-bound.  Another example where this behaviour might arise is as follows.  Suppose we seek a minimiser to a variational problem whose Euler-Lagrange,  or other variational equations are sufficiently degenerate that known methods do not provide existence or regularity.  Suppose further one can perturb these equations to gain ellipticity and thereby existence and some regularity.  The solutions to the perturbed equation may have a weakly convergent subsequence for which local uniform bounds on nonlinear quantities hold and which do not depend on the ellipticity constants.  Strong convergence may then imply that a weak limit satisfies the unperturbed equation from which one might deduce higher regularity.  Despite all the suppositions here,  we outline some concrete examples later among  mean distortion functionals and Ahlfors-Hopf type equations.

\subsection{Finite distortion functions and polyconvexity}
Let $\Omega$ be a planar domain. A mapping $f:\Omega\to\IC$ has finite distortion if 
\begin{enumerate}
\item $f\in W^{1,1}_{loc}(\Omega)$,  the Sobolev space of functions with locally integrable first derivatives,
\item the Jacobian determinant $J(z,f)\in L^{1}_{loc}(\Omega)$, and 
\item there is a measurable function $ K(z,f)\geq 1$, finite almost everywhere, such that 
 \begin{equation}\label{1.1}
 |Df(z)|^2 \leq K(z,f) \, J(z,f), \hskip10pt \mbox{ almost everywhere in $\Omega$}.
 \end{equation}
\end{enumerate}
We recommend \cite[Chapter 20]{AIM} for the basic theory of mappings of finite distortion and the associated governing equations; degenerate elliptic Beltrami systems.  

\medskip

In (\ref{1.1}) the operator norm of $Df$ is used,  however for variational problems it is more common to use 
\begin{equation}\label{IK} \IK(z,f)=\frac{\|Df(z)\|^2}{J(z,f)}, \end{equation}
 where $\|\cdot\|$ is the Hilbert-Schmidt norm.

The notion of polyconvexity was introduced to the theory on nonlinear elasticity by Ball in \cite{Ball},  and has proved an important concept in the calculus of variations ever since.  A matrix function $\Xi:\IR^{n\times n}\to\IR$ is polyconvex if it can be expressed as a convex function of minors of $\IR^{n\times n}$. See \cite[Definition 10.25]{MR} and in the current context \cite[Section 14.2]{AIMO}. Here we only consider the two-dimensional case and our main interest here is the convexity of the function $\frac{x^2}{y}$ (which the reader should compare with the definition at (\ref{IK})). Namely,
\[
\frac{x^2}{y}-\frac{x_0^2}{y_0}\geq\frac{2x_0}{y_0}(x-x_0)-\frac{x_0^2}{y_0^2}(y-y_0),
\]
for any $x,x_0\geq0$ and $y,y_0>0$. In our applications $x$ will be $|Df(z)|$ or $\|Df(z)\|$, and $y$ will be $J(z,f)$. We state our first lemma in a more general setting, but the proof is essentially the same as \cite[Theorem 12.2]{AIMO}. Also we only write for $|Df(z)|$, but everything follows similarly for $\|Df(z)\|$.

\begin{lemma}\label{polyconvexity} Let $\Omega\subset\IC$, and let 
\begin{itemize}  
\item $\{f_j\}_{j=1}^{\infty}$ be a sequence of $W_{loc}^{1,q}(\Omega)$ functions, where $1\leq q<\infty$, 
\item $f$ be a weak limit of $f_j$ in $W_{loc}^{1,q}(\Omega)$, 
\item $J(z,f_j)\rightharpoonup J(z,f)$ weakly in $L_{loc}^1(\Omega)$, 
\item $J(z,f)>0$ almost everywhere in $\Omega$, and 
\item $\Phi(x,y):\IR^+\cup\{0\}\times\IR^+\cup\{0\}\to\IR^+\cup\{0\}$ be a convex function which has partial derivatives $\Phi_x$, $\Phi_y$ almost everywhere, and
\[
0\leq\Phi_x(|Df(z)|,J(z,f))<\infty,\quad\Big|\Phi_y(|Df(z)|,J(z,f))\Big|<\infty,
\]
for almost every $z\in\Omega$.
\end{itemize} 
Then
\[
\int_\Omega\Phi(|Df(z)|,J(z,f))\leq\liminf_{j\to\infty}\int_\Omega\Phi(|Df_j(z)|,J(z,f_j)).
\]
\end{lemma}
In \cite[Theorem 12.2]{AIMO}  equality is actually proved for certain functionals  because $f_j$ is assumed a minimising sequence. This should remind us of the Radon-Riesz property. 

\subsection{Radon-Riesz Property.}
A Banach space is called a Radon-Riesz space, if every weakly convergent sequence $x_j\rightharpoonup x$ with $\|x_j\|\to\|x\|$ is strongly convergent, namely $\|x_j-x\|\to0$. See e.g. \cite{Megginson}.
\begin{lemma}\label{RadonRiesz}
Every uniformly convex Banach space is a Radon-Riesz space. In particular, every $L^p$ space  with $1<p<\infty$ is a Radon-Riesz space.
\end{lemma}
Clearly the map 
\[ f\mapsto\int_\Omega\Phi(|Df(z)|,J(z,f))\; dz \] is usually not a norm. Our aim is to prove that under certain criteria, weak convergence of $f_j\to f$ implies convergence strongly in some $W_{loc}^{1,q}(\Omega)$. This is our main result.

\begin{theorem}\label{Maintheorem}
Let $\Omega\subset\IC$ be a domain.  Suppose that
\begin{itemize}
\item $f_j$ is a sequence of $W_{loc}^{1,q}(\Omega)$ functions, for some $1<q<\infty$, 
\item $f$ is a weak limit of $f_j$ in $W_{loc}^{1,q}(\Omega)$, 
\item $J(z,f_j)\rightharpoonup J(z,f)$ weakly in $L_{loc}^1(\Omega)$, 
\item $J(z,f)>0$ almost everywhere  in $\Omega$, and 
\item $\Phi$, $\Phi_j:\IR^+\cup\{0\}\times\IR^+\cup\{0\}\to\IR^+\cup\{0\}$ is a sequence functions which satisfies the following conditions:
\begin{enumerate}
\item There is a $p>1$ such that $\Phi_j(|Df_j|,J(z,f_j))$ are uniformly bounded in $L^p(\Omega)$, and
\begin{equation}\label{PhiConv}
\lim_{j\to\infty}\int_\Omega\Phi_j^p(|Df_j(z)|,J(z,f_j))=\int_\Omega\Phi^p(|Df(z)|,J(z,f))
\end{equation}
\item For almost every pair of $x,y$, $\Phi_j(x,y)$ is a non-decreasing sequence and $\Phi_j(x,y)\to\Phi(x,y)$.
\item Each $\Phi_j(x,y)$ is a convex function and has partial derivatives for almost every pair $x,y$, and
\[
0\leq(\Phi_j)_x(|Df(z)|,J(z,f))<\infty,\quad\Big|(\Phi_j)_y(|Df(z)|,J(z,f))\Big|<\infty,
\]
for almost every $z\in\Omega$.
\item There is an $s\in(0,1-\frac{1}{p})$ such that each $\Phi_j(x,y)y^s$ is a convex function.
\item If $\Phi_j(|Df_j(z)|,J(z,f_j))\to\Phi(|Df(z)|,J(z,f))$ and $J(z,f_j)\to J(z,f)$ pointwise almost everywhere, then $|Df_j(z)|\to|Df(z)|$ pointwise almost everywhere.
\end{enumerate}
\end{itemize} 
Then, there is a subsequence for which the following convergence holds both strongly and pointwise almost everywhere.
\begin{itemize}
\item $\Phi_j(|Df_j|,J(z,f_j))\to\Phi(|Df|,J(z,f))\mbox{ in }L^p(\Omega)$,
\item $f_j\to f\mbox{ in }W_{loc}^{1,r}(\Omega),\quad0<r<q$,
\item $J(z,f_j)\to J(z,f)\mbox{ in }L_{loc}^r(\Omega),\quad0<r<1$,
\item $\mu_{f_j}\to\mu_f\mbox{ in }L_{loc}^r(\Omega),\quad0<r<\infty$,
\end{itemize}
where $\mu_f=f_\zbar/f_z$ is the Beltrami coefficient of $f$.
\end{theorem}

A weighted case can be proved in a similar way:
\begin{theorem}\label{TheoremWeighted}
Let $f$, $f_j$, $\Phi$, $\Phi_j$ be same as in Theorem 1, $\eta_j, \eta>0$ a.e. in $\Omega$, both in $L_{\loc}^\infty(\Omega)$, $\eta_j\to\eta$ locally uniformly in $\Omega$, and (\ref{PhiConv}) be modified to
\begin{equation}\label{ConvWeighted}
\lim_{j\to\infty}\int_\Omega\Phi_j^p(|Df_j(z)|,J(z,f_j))\eta_j(z)=\int_\Omega\Phi^p(|Df(z)|,J(z,f))\eta(z)
\end{equation}
Then the same results in Theorem 1 hold.
\end{theorem}
\section{Proof of Theorem \ref{Maintheorem} and \ref{TheoremWeighted}.}
To simplify notation we write
\begin{align*}
\Phi_{k,j}&=\Phi_k(|Df_j(z)|,J(z,f_j)), & \Phi_{k,f}&=\Phi_k(|Df(z)|,J(z,f)), & \Phi_f&=\Phi(|Df(z)|,J(z,f)) \\
 J_j& =J(z,f_j), & J_f&=J(z,f).
\end{align*}
\begin{lemma}\label{lemma3}
$\Phi_{j,j}\rightarrow\Phi_f$ strongly in $L^p(\Omega)$.
\end{lemma}
\noindent{\bf Proof.} By (\ref{PhiConv}), the sequence $\Phi_{j,j}$ is uniformly bounded in $L^p(\Omega)$, so there is a weak limit $\Psi\in L^p(\Omega)$. Note Lemma \ref{polyconvexity} holds for every fixed $k$ and any measurable subset $\Omega'\subset\Omega$. Thus
\[
\int_{\Omega'}\Phi_{k,f}\leq\liminf_{j\to\infty}\int_{\Omega'}\Phi_{k,j}\leq\liminf_{j\to\infty}\int_{\Omega'}\Phi_{j,j}=\int_{\Omega'}\Psi.
\]
This holds for every $k$, and so by Fatou's lemma,
\[
\int_{\Omega'}\Phi_f\leq\liminf_{k\to\infty}\int_{\Omega'}\Phi_{k,f}\leq\int_{\Omega'}\Psi.
\]
From the  Lebesgue differentiation theorem, $\Phi_f(z)\leq\Psi(z)$ for almost every  $z\in\Omega$. Now it follows from (\ref{PhiConv}) that
\[
\int_\Omega\Phi_f^p\leq\int_\Omega\Psi^p\leq\liminf_{j\to\infty}\int_\Omega\Phi_{j,j}^p=\int_\Omega\Phi_f^p.
\]
This forces $\Phi_f=\Psi$, so it is the weak limit of $\Phi_{j,j}$ in $L^p(\Omega)$, and then the claim follows from Lemma \ref{RadonRiesz}.\hfill $\Box$

\bigskip

Next, as $J_f>0$ almost everywhere we can choose $\Omega_\varepsilon\Subset\Omega$ on which
\[
\varepsilon<J_f<\frac{1}{\varepsilon},\quad\Phi_f<\frac{1}{\varepsilon},
\]
and
\[
\Big|\Omega-\bigcup_{\varepsilon>0}\Omega_\varepsilon\Big|=0.
\]
Now let $s\in(0,1-\frac{1}{p})$ be as in Condition 4. We choose a $p'\in(1,p)$ such that $sp'<(s+\frac{1}{p})p'<1$. Then
\begin{lemma}\label{lemma4}
\[
\lim_{j\rightarrow\infty}\int_{\Omega_\varepsilon}\Phi_{j,j}^{p'}J_j^{sp'}=\int_{\Omega_\varepsilon}\Phi_f^{p'}J_f^{sp'}.
\]
\end{lemma}
\noindent{\bf Proof.} One direction comes from the polyconvexity which is same as above:
\[
\int_{\Omega_\varepsilon}\Phi_f^{p'}J_f^{sp'}\leq\liminf_{k\to\infty}\int_{\Omega_\varepsilon}\Phi_{k,f}^{p'}J_f^{sp'}\leq\liminf_{j\rightarrow\infty}\int_{\Omega_\varepsilon}\Phi_{j,j}^{p'}J_j^{sp'}.
\]
Lemma \ref{lemma3} and the choice of $p'$ gives
\[
\lim_{j\rightarrow\infty}\int_{\Omega_\varepsilon}\big|\Phi_{j,j}^{p'}J_j^{sp'}-\Phi_f^{p'}J_j^{sp'}\big|\leq\lim_{j\rightarrow\infty}C\|\Phi_{j,j}-\Phi_f\|_{L^p(\Omega_\epsilon)}^{p'}\|J_j\|_{L^1(\Omega_\epsilon)}^{sp'}=0.
\]
So we only need to show
\begin{equation}\label{equation3}
\limsup_{j\rightarrow\infty}\int_{\Omega_\varepsilon}\Phi_f^{p'}J_j^{sp'}\leq\int_{\Omega_\varepsilon}\Phi_f^{p'}J_f^{sp'}.
\end{equation}
Note that as $sp'<1$, the function $x\mapsto x^{sp'}$ is concave for $x>0$. So
\[
J_j^{sp'}-J_f^{sp'}\leq sp'J_f^{sp'-1}(J_j-J_f).
\]
It follows that
\[
\int_{\Omega_\varepsilon}\Phi_f^{p'}(J_j^{sp'}-J_f^{sp'})\leq sp'\int_{\Omega_\varepsilon}\Phi_f^{p'}J_f^{sp'-1}(J_j-J_f)\rightarrow0,
\]
as $J_j\rightharpoonup J_f$ in $L^1(\Omega)$. This proves (\ref{equation3}) and completes the proof of the lemma.\hfill $\Box$

\bigskip

If we apply Lemma \ref{lemma4} and follow the same proof as in Lemma \ref{lemma3} we can establish the following lemma.
\begin{lemma}\label{lemma5}
$\Phi_{j,j}J_j^s\to\Phi_fJ_f^s$ strongly in $L^{p'}(\Omega_\varepsilon)$
\end{lemma}

\bigskip

Now by Lemma \ref{lemma5}, in every $\Omega_\varepsilon$, up to a subsequence we have the pointwise almost everywhere convergence $\Phi_{j,j}J_j^s\to\Phi_f J_f^s$. Let $\varepsilon\to0$, we can choose diagonally and to obtain a subsequence that converges pointwise in $\Omega$. By Lemma \ref{lemma3} we also have the pointwise convergence $\Phi_{j,j}\to\Phi_f$. It follows that $J_j\to J_f$ pointwise almost everywhere, then $|Df_j|\to|Df|$ pointwise almost everywhere, and then $K(z,f_j)\to K(z,f)$, $|\mu(z,f_j)|\to|\mu(z,f)|$, $|(f_j)_z|\to|f_z|$ and $|(f_j)_\zbar|\to|f_\zbar|$, all pointwise almost everywhere.

Next, let $\Omega'\Subset\Omega$. By Vitali's convergence theorem, both $|(f_j)_z|\to|f_z|$ and $|(f_j)_\zbar|\to|f_\zbar|$ strongly in $L^r(\Omega')$, for any $1\leq r<q$. In particular,
\[
\int_{\Omega'}|(f_j)_z|^r\to\int_{\Omega'}|f_z|^r,\quad\int_{\Omega'}|(f_j)_\zbar|^r\to\int_{\Omega'}|f_\zbar|^r.
\]
On the other hand, $Df_j\rightharpoonup Df$ weakly in $L^r(\Omega')$. So it follows from Lemma \ref{RadonRiesz} that $Df_j\to Df$ strongly in $L^r(\Omega')$.  We may now apply Vitali's convergence theorem once again to get the remaining claims of Theorem \ref{Maintheorem}.\\

The proof of Theorem \ref{TheoremWeighted} is essentially same. We only need to replace $\Phi$ by $\Phi\eta^\frac{1}{p}$ in the proof.

\section{Applications}
\subsection{$\exp(p)$ minimising sequence}
In \cite[Theorem 12.2]{AIMO} it is proved that the $\exp(p)$ mean distortion
\[ f \mapsto \int_\Omega\exp[p\IK(z,f)] \; dz\]
for homeomorphisms from $\overline{\Omega}$ to $\overline{\Omega'}$ admits a minimiser for suitable boundary data. In fact, there is  minimising sequence $f_j$ which converges weakly to a minimiser $f$ in the Sobolev-Orlicz space $W^{1,P}(\Omega)$, where $P(t)=\frac{t^2}{\log(e+t)}$, and
\[
\int_\Omega\exp[p\IK(z,f)]=\lim_{j\to\infty}\int_\Omega\exp[p\IK(z,f_j)].
\]
Thus in Theorem \ref{Maintheorem} we may set $\Phi_j(x,y)=\Phi(x,y)=\exp[\frac{p}{2}\frac{x^2}{y}]$ and any small $s>0$, to obtain the strong convergence of $\exp(p\IK(z,f_j))$, $Df_j$, and so forth. In particular, for the case $p>2$ we have $f_j\to f$ strongly in $W_{loc}^{1,2}(\Omega)$ (cf. \cite{AGRS}). Furthermore, after changing variables we also have
\[
\int_{\Omega'}\exp[p\IK(w,h)]J(w,h)=\lim_{j\to\infty}\int_{\Omega'}\exp[p\IK(w,h_j)]J(w,h_j),
\]
where $h_j=f_j^{-1}$, $h=f^{-1}$. Again in Theorem \ref{Maintheorem} we may set $\Phi_j(x,y)=\Phi(x,y)=\exp[\frac{p}{2}\frac{x^2}{y}]y^\frac{1}{2}$, to obtain the strong convergence of $h_j\to h$ in $W_{loc}^{1,2}(\Omega')$.

\subsection{$L^p$ minimising sequence}
In \cite {MY1} we considered the boundary value problems for $L^p$- mean distortion for self homeomorphisms of the unit disk $\ID$,
\[ f \mapsto \int_\ID \IK^p(z,f)\; dz. \] 
Unfortunately, in this case $Df$ is not apriori in a sufficiently regular space, so a minimising sequence $f_j$ might not have weakly convergent $J(z,f_j)$ in $L_{loc}^1(\ID)$, and further the limit function $f$ might not be a homeomorphism. However, the inverse sequence $h_j\rightharpoonup h$ weakly in $W^{1,2}(\ID)$, where $h$ is a minimiser in an enlarged space where pseudo-inverses exist. In particular,
\[
\int_\ID\IK^p(w,h)J(w,h)=\lim_{j\to\infty}\int_\ID\IK^p(w,h_j)J(w,h_j).
\]
Here we can apply Theorem \ref{Maintheorem}. Set $\Phi_j(x,y)=\Phi(x,y)=(\frac{x^2}{y})y^\frac{1}{p}$ and any small $s>0$ to obtain the strong convergence of $h_j\to h$ in $W^{1,2}(\ID)$. In fact more is true. In the enlarged space there are pseudo-inverses $f_j=h_j^{-1}$ and $f=h^{-1}$, $f_j\rightharpoonup f$ in $W^{1,\frac{2p}{p+1}}(\ID)$, each $f_j(\ID)\subset\ID$, and
\[
\int_\ID\IK^p(z,f)=\lim_{j\to\infty}\int_\ID\IK^p(z,f_j).
\]
We consider the sequence $J^{\frac{1}{2}}(z,f_j)$, which is bounded in $L^2(\ID)$. For all $\varphi\in C_0^\infty(\ID)$,
\begin{align*}
\lim_{j\to\infty}\int_\ID J^{\frac{1}{2}}(z,f_j)\varphi(z)&=\lim_{j\to\infty}\int_\ID J^{\frac{1}{2}}(w,h_j)\varphi(h_j(w))\\
&=\int_\ID J^{\frac{1}{2}}(w,h)\varphi(h(w))=\int_\ID J^{\frac{1}{2}}(z,f)\varphi(z).
\end{align*}
So $J^{\frac{1}{2}}(z,f)$ is the weak limit of $J^{\frac{1}{2}}(z,f_j)$ in $L^2(\ID)$. Also the minimiser $h$ has the holomorphic Hopf differential
\[
\Phi=\IK^{p-1}(w,h)h_w\overline{h_\wbar}=\IK^p(w,h)J(w,h)\frac{\overline{\mu(w,h)}}{1+|\mu(w,h)|^2}.
\]
This gives us two cases: firstly, $\Phi\equiv0$, then $h$ is a holomorphic function. But we know $h$ is monotone, so actually it is a conformal mapping, thus $f(\ID)=\ID$. If $\Phi$ is not identically zero,  then  $J(w,h)>0$ almost everywhere in $\ID$. But we know $J(w,h)=0$ almost everywhere in $\ID-f(\ID)$. So in either case we have $|f(\ID)|=\pi$. Now
\[
\pi\geq\liminf_{j\to\infty}\int_\ID J(z,f_j)\geq\int_\ID J(z,f)=\int_{f(\ID)}\frac{1}{J(w,h)}J(w,h)=\pi.
\] 
So the inequalities hold with equality and again Radon-Riesz Lemma \ref{RadonRiesz} implies the strong convergence of $J(z,f_j)\to J(z,f)$ in $L^1(\ID)$. Now Theorem \ref{Maintheorem} applies and we also get that $\IK(z,f_j)\to\IK(z,f)$ strongly in $L^p(\ID)$, and $f_j\to f$ strongly in $W^{1,q}(\ID)$ for all $1\leq q<\frac{2p}{p+1}$.

\subsection{Truncated exponential minimising sequence}
The exponential finite distortion problem is not variational \cite{MY3}.  Thus in \cite {MY2}, to study the exponential finite distortion problems,  we consider the truncated problems and the associated inverse problems:
\[
f\mapsto \int_\ID\sum_{n=0}^N\frac{p^n\IK^n(z,f)}{n!},\quad h\mapsto \int_\ID\sum_{n=0}^N\frac{p^n\IK^n(w,h)}{n!}J(w,h).
\]
As linear combinations of the $L^p$ problems, in the enlarged space there are minimisers $h_N$ which have holomorphic Ahlfors-Hopf differentials
\[
\Psi_N=\sum_{n=0}^N\frac{p^n\IK^n(w,h_N)}{n!}(h_N)_w\overline{(h_N)_\wbar}.
\]
In Theorem $\ref{Maintheorem}$ we set
\[
\Phi_{N,N}=\sqrt{\sum_{n=0}^N\frac{p^n\IK^n(w,h_N)}{n!}J(w,h_N)}
\]
which is a polyconvex function of $Dh_N$ for sufficiently large $N$. So Theorem $\ref{Maintheorem}$ gives $h_N\to h$ strongly in $W^{1,2}(\ID)$. A similar argument as in the last part tells us that their inverses $f_N=h_N^{-1}$ also converge strongly to $f=h^{-1}$, but in $W^{1,q}(\ID)$ for all $q\in[1,2)$. In fact we can prove the limit function $f$ is a homeomorphic minimiser for the $\exp(p)$ problem, and if the sequence $\frac{\Psi_N}{\|\Psi_N\|_{L^1(\ID)}}$ is nondegenerate, then there is a holomorphic $\Psi$ such that $\Psi_N\to\Psi$ locally uniformly, and then by the strong convergence we have
\[
\Psi=\lim_{N\to\infty}\Psi_N=\lim_{N\to\infty}\sum_{n=0}^N\frac{p^n\IK^n(w,h_N)}{n!}(h_N)_w\overline{(h_N)_\wbar}=\exp[p\IK(w,h)]h_w\overline{h_\wbar},
\]
and we can prove such an $h$ is diffeomorphic in $\ID$, and then so is $f=h^{-1}$. See \cite {MY2} for detailed discussions.\\

\section{Mappings between surfaces.}

In this section we consider the following minimising problem. Let $f_0:S\to\tilde{S}$ be a quasiconformal mapping,  where $S$, $\tilde{S}$ are compact Riemann surfaces.  We seek to study the critical points of 
\begin{equation}\label{surface} f \mapsto \int_S\exp[p\IK(z,f)]\; d\sigma(z) \end{equation}
where $f$ is in the same homotopy class as the datum $f_0$ and where $d\sigma(z)$ is the hyperbolic area measure on $S$. 

\medskip

Again,  this problem is not variational.  However,  following the argument of the previous section we can show that the inverse of a minimiser satisfies the variational equations.  That is the Ahlfors-Hopf equation.  

We lift this problem to the universal cover.  Any such mapping $f$ as appears in (\ref{surface}) has a lift $\tilde{f}:\ID\to\ID$ which commutes with the fundamental groups $\Gamma$, of $S$, and $\tilde{\Gamma}$, of $\tilde{S}$.  That is ${f_0}_*:\pi_1(S)\to \pi_1(\tilde{S})$ induces an isomorphism which we simply denote as $\gamma\mapsto\tilde{\gamma}$.  Then
\begin{equation}\label{7}
\tilde{f} \circ \gamma = \tilde{\gamma}\circ \tilde{f} :\ID\to\ID.
\end{equation} 
Under these circumstances $\tilde{f}$ extends to $\IS=\partial\ID$ and $\tilde{f}|\IS$ is quasisymmetric.  Notice too that if $\tilde{f_j}\to \tilde{f}$ locally uniformly,  and if $\tilde{f_j}$ satisfies (\ref{7}), then so does the limit $\tilde{f}$.

\medskip

Next, if $\mathcal{P}$ denotes a (hyperbolically) convex fundamental polyhedron for $\Gamma$,  then
\begin{equation}  \int_S\exp[p\IK(z,f)]\; d\sigma(z) = \int_\mathcal{P} \exp[p\IK(z,\tilde{f})] \eta(z) \; dz
\end{equation}
where the weight $\eta(z)$ is  the hyperbolic metric $\frac{1}{(1-|z|^2)^2}$. See \cite{Katok} for more details. 

Thus we may now we can consider (writing $f$ for $\tilde{f}$) 
\[
\int_\mathcal{P}\exp[p\IK(z,f)]\eta(z),
\]
where  we also impose the automorphic condition at (\ref{7}). Again we can consider the inverse problem
\[
\int_{\tilde{\mathcal{P}}}\exp[p\IK(w,h)]J(w,h)\eta(h),
\]
where $\tilde{\mathcal{P}}$ is a fundamental domain for $\tilde{\Gamma}$ and we impose the automorphy condition 
\begin{equation}\label{8}
\gamma \circ h= h\circ \tilde{\gamma} :\ID\to\ID.
\end{equation} 
The associated truncated problems are to minimise
\begin{equation}\label{9}
\int_{\tilde{\mathcal{P}}}\sum_{n=0}^N\frac{p^n\IK^n(w,h)}{n!}J(w,h)\eta(h).
\end{equation} 
Then, as we have found earlier,  there is a sequence of minimisers $h_N$ (in the enlarged space) with holomorphic Ahlfors-Hopf differentials
\begin{equation}\label{NHopf}
\Psi_N=\sum_{n=0}^N\frac{p^n\IK^n(w,h_N)}{n!}(h_N)_w\overline{(h_N)_\wbar}\eta(h_N).
\end{equation}
In fact this is basically how Ahlfors sets up his approach to the proof of Teichm\"uller's theorem.  However he multiplies through by a ``convergence factor" to eliminate the bad term $\eta(h)$ in (\ref{9}) as it  will make no difference in his application as he lets $p\to\infty$.  Now $\eta(h_N)\to\eta(h)$ locally uniformly in $\ID$ and so uniformly on $\tilde{\mathcal{P}}$. So by Theorem \ref{TheoremWeighted} there is an $h$ such that $h_N\to h$ strongly in $W^{1,2}(\tilde{\mathcal{P}})$. Here the difference is that we know the space of quadratic differentials is finite dimensional by the Riemann-Roch theorem, so $\frac{\Psi_N}{\|\Psi_N\|_{L^1(\mathcal{P})}}$ is nondegenerate, and so there is another holomorphic $\Psi$ such that
\[
\Psi=\exp[p\IK(w,h)]h_w\overline{h_\wbar}\eta(h).
\]
In fact this argument works for each domain $\tilde{\gamma}(\tilde{\mathcal{P}})$, where $\tilde{\gamma}\in\tilde{\Gamma}$, so both of the functions $h$ and $\Psi$ extend to $\ID$, where $\Psi$ is a holomorphic function in $\ID$, and the equation $\Psi=\exp[p\IK(w,h)]h_w\overline{h_\wbar}\eta(h)$ holds over $\ID$. 

Also, since $f$ is in the exponential class and is automorphic with respect to Fuchsian groups of compact type,  both $f$ and $h$ must be self-homeomorphisms of the closed disk $\overline{\ID}$.  To see this we may argue as follows. If $\ID=h(\Omega)$,  $h:\Omega\to \ID$ smooth,  then
$h(\Omega) = (\tilde{\gamma}\circ h \circ \gamma)(\Omega)$,  
$ \ID=\tilde{\gamma}^{-1}(\ID) =  h(\gamma(\Omega))$.
As $h$ is a diffeomorphism onto,  it follows that $\gamma(\Omega)=\Omega$ for each $\gamma\in \Gamma$.  As $\Omega$ is simply connected $\Omega/\Gamma$ is a Riemann surface with fundamental group $\Gamma$. This is the fundamental group of a closed surface,  so $\Omega/\Gamma$ is a closed surface and a fundamental domain for $\Gamma$ lies in $\Omega$.  This shows $\Omega=\ID$.

As noted,  the automorphy condition passes to the limit.  In the case $\eta(z)=\frac{1}{(1-|z|^2)^2}$, this and (\ref{NHopf}) together gives
\begin{align*}
\Psi_N(w)&=\sum_{n=0}^N\frac{p^n\IK^n(w,\gamma\circ h_N)}{n!}(\gamma\circ h_N)_w\overline{(\gamma\circ h_N)_\wbar}\eta(\gamma\circ h_N)\\
&=\sum_{n=0}^N\frac{p^n\IK^n(w,h_N\circ\tilde{\gamma})}{n!}(h_N\circ\tilde{\gamma})_w\overline{(h_N\circ\tilde{\gamma})_\wbar}\eta(h_N\circ\tilde{\gamma})=\Psi_N(\tilde{\gamma})\tilde{\gamma}'^2.
\end{align*}
By the strong convergence these properties of $h_N$ and $\Psi_N$ persist in the limit for $h$ and $\Psi$. We conclude as follows:
\begin{theorem}
Consider the inverse hyperbolic exponential finite distortion problem $h\mapsto\int_\ID\exp[p\IK(w,h)]J(w,h)\eta(h)$, where $\eta(z)=\frac{1}{(1-|z|^2)^2}$, and $h$ is a self-homeomorphism of $\overline{\ID}$ and is automorphic with respect to the Fuchsian groups $(\tilde{\Gamma},\Gamma)$. Then, there is a critical point $h$ and holomorphic Ahlfors-Hopf differential $\Psi$ such that
\begin{equation}\label{12}
\Psi=\exp[p\IK(w,h)]h_w\overline{h_\wbar}\eta(h).
\end{equation}
Furthermore, $\Psi=\Psi(\tilde{\gamma})\tilde{\gamma}'^2$ for every $\tilde{\gamma}\in\tilde{\Gamma}$.
\end{theorem}
Note that in these circumstances,  we prove in \cite{MY2} that there is a diffeomorphic solution to (\ref{9}).  Suitably normalised,  we expect this solution to be unique.  Further,  we also show that any quasiconformal solution to (\ref{9}) is already a diffeomorphism.  These two fact strongly suggest the minimiser is a diffeomorphism,  but our results to date fall a little short of this as our uniqueness results are not yet strong enough.

\end{document}